\documentclass[11pt,oneside,reqno]{amsart}

\usepackage{amsmath, amsrefs}
\usepackage{amsfonts}
\usepackage{amssymb}
\usepackage{amsthm}
\usepackage{bbm}
\usepackage{color, comment}
\usepackage{enumitem}
\usepackage[a4paper,textwidth=14cm]{geometry}
\usepackage[hidelinks]{hyperref}
\usepackage{tikz}
\usetikzlibrary{arrows,shadows, shapes, decorations, decorations.markings}
\usepackage[a4paper,textwidth=14cm]{geometry}
\usepackage[nameinlink]{cleveref}
\hypersetup{
    colorlinks=true,
    linkcolor=magenta,
    filecolor=magenta,
}

\newcommand{\RZ}{\mathbb{R}/\mathbb{Z}}
\newcommand{\eps}{\varepsilon}
\renewcommand{\H}{\mathcal{H}}

\newcommand{\R}{\mathbb{R}}
\newcommand*{\genbf}[1]{\ifmmode\mathbf{#1}\else\textbf{#1}\fi}

\newtheorem{theorem}{Theorem}[section]
\newtheorem{lemma}[theorem]{Lemma}
\newtheorem{proposition}[theorem]{Proposition}
\newtheorem{corollary}[theorem]{Corollary}

\theoremstyle{definition}
\newtheorem{definition}[theorem]{Definition}
\newtheorem{remark}[theorem]{Remark}
\newtheorem{example}[theorem]{Example}
\newtheorem{problem}{Problem}

\title{The planar Plateau's problem via capillarity}
\author{Kennedy Obinna Idu}
\address{Department of Mathematics, University of Toronto, 
			Toronto, Ontario M5S 2E4, Canada.}
\email{o.idu@utoronto.ca}
\date{}

\begin{document}
\begingroup
\def\uppercasenonmath#1{} % this disables uppercasing title
\let\MakeUppercase\relax % this disables uppercasing authors
\maketitle
\vspace{-1.5em}
\endgroup

\begin{abstract} 
The Plateau's problem seeks to determine a surface of minimal area which spans a given boundary. It is widely studied for its varied mathematical formulations, applications and relevance to physical models such as soap films. We revisit the problem and study a soap film model in the spirit of capillarity formulations in two dimensions. Our approach introduces a nonlocal geometric potential in the variational length minimization scheme. This incorporates effects of thickness of soap films and provides insight into addressing the so-called collapsing phenomenon and other observable physical phenomena and properties. 
\end{abstract}
% \textbf{MSC (2010):} 28A75, 28A78, 26A16.\\
% \textbf{Keywords:} ....
%Note: you can toggle the labels by commenting 
%\begin{verbatim}\usepackage{showkeys}\end{verbatim}

%\tableofcontents

\section{Introduction}
%Most spatial arrangements in nature tend to capture the principal focus of variational problems -- \textit{area} minimisation subject to topological constraints.
%The most widely studied is, perhaps, that of soap films. However, we have an avalanche of models from the modeling of repulsive forces of fibres, charged loop of string and of macromolecules such as DNA and proteins which seek to attain a minimum state of a suitable energy. 

%As in soap films, not only steric arrangements but also, rheological factors play a vital role in the minimisation of interfacial cell areas of most spatial structures. It is often essential to establish embeddability of these structures and to guarantee control of unwanted self-intersections in their accompanying variational analysis.

The famous Plateau problem entails finding a \textit{surface} of minimal \textit{area} which \textit{spans} a given \textit{boundary}. It admits quite a number of formulations (precipitating from the various notions of a surface, area, and what it means to span a boundary) which attempt to describe physical properties of soap films. A more precise formulation is: \textit{Let $\Gamma\subset \R^{n+1}$ be a $(n-1)$-dimensional surface  without boundary. Find a $n$-dimensional surface $\Sigma$ such that:}
\begin{equation}\label{Plateau1}
    H_\Sigma=0, \qquad \partial \Sigma=\Gamma,
\end{equation}
where $H_\Sigma$ is the mean curvature of $\Sigma$ and $\partial \Sigma$ represents the boundary, and by \textit{find} we  mean to ``prove the existence of".

In \cite{GD14}, the author raised striking inquiries into trends in advances on the problem. A chief issue was recovering the original model in the sense of the motivating physical phenomena of soap films (see also \cites{Har04, HP16Spr}). The difficulty has been, in part, due to finding the right notion for what it means for a surface to ``span a given boundary''. Recent progress have been set forth by the seminal works of Harrison and Pugh \cites{HP16,HP17}. In \cite{HP16}, the authors define a topological notion of \textit{span} using a homotopy class condition. This condition in the Harrison-Pugh formulation of the Plateau problem proves to be very crucial in capturing the notion of span of boundaries which resemble soap films by \textit{modding out} possible ``droplets'' known to occur at points of maximal curvature on the boundary under natural spanning assumptions in the so-called capillarity formulation of the problem (see e.g., \cite{MM16}). Using the Harrison-Pugh model, King \textit{et al.} \cite{KMS20} insightfully confront another hitch in the mathematical theory of recovering the physical soap film model in the Plateau problem: \textit{length scale} property. Indeed, the formulation in \eqref{Plateau1} fails to capture the physical observation of the disappearance of a soap film by arbitrarily scaling the boundary wire. King \textit{et al.} borrow from the insight of Maggi \textit{et al.} \cite{MSS19} on the study of \textit{almost-minimal surfaces} and consider the modelling of soap films as regions $E\subset \R^{n+1}$ with small but positive volume $|E|=\eps$. This is the  framework of the Capillarity problem and tacitly factors a length scale in the Plateau problem which is recovered in the vanishing volume limit.

\indent In the planar setting, the Plateau problem admits the following formulation:
\begin{problem}\label{Prob:Plateau}
Let $S\subset\R^2$ be a finite set. Find a compact and connected set $K$ of minimal $\H^1$-measure that contains $S$.
\end{problem}

This coincides with the classical Euclidean Steiner tree problem (we refer to \cite{Paolini12Steiner} and the references therein).

In this paper, we are interested in a solution to \Cref{Prob:Plateau}. We consider a different approach in the spirit of the Capillarity model as in Maggi \textit{et al.} \cites{MSS19,KMS20}. We note that in dimension two, the capillarity model involves curves with enclosed area.
Our insight is to model thickness by penalising self-intersection through the introduction of a nonlocal geometric potential when minimising length in the class of rectifiable curves that enclose a fixed area.\\ The precise formulation is as follows:
%Our approach entails the minimisation problem \eqref{MP} with respect to the functional $G_\delta$ consisting of a nonlocal energy penalising self-intersection in the class of rectifiable curves enclosing a fixed area.

\subsection{Variational formulation of problem}
\indent We fix $\eps>0$ and for any positive integer $N\ge 2, \{x_1,\dots,x_N\}\subset \R^2$, and consider the class $\mathcal{E}$ of all closed Lipschitz curves $\gamma: \RZ \rightarrow \mathbb{R}^2$ satisfying
\begin{enumerate}
\item[E1.] Inclusion of finite set: $\gamma(I) \supset \{x_1,\ldots, x_N\}$, 
\item[E2.] Well-defined interior and exterior: Winding number (see \Cref{Def: WindingNum}) $\omega(\gamma, x)\in \{0,1\}$ for any $x\notin \gamma(\RZ)$,
\item[E3.] Positive area enclosure:  $\left|\mathtt{int}(\gamma):=\{x\in\R^2: \omega(\gamma,x)=1\}\right|=\eps$, 
\item[E4.] $\gamma$ has constant speed.
\end{enumerate}

For any $\delta>0$, $s\in (0,1)$ and $1\le p\le \infty$, we are interested in
\begin{subequations}
 \begin{equation}\label{Main}
\min\limits_{\gamma\in \mathcal{E}} \left\{G_\delta(\gamma):= \ell(\gamma) + \delta \int_{\RZ}\int_{\RZ}\frac{\left|\tau(t)-\tau(t')\right|^p}{|\gamma(t)-\gamma(t')|^{1+sp}}|\dot{\gamma}(t)||\dot{\gamma}(t')|dtdt'\right\},
\end{equation}
where $\tau(t):=\frac{\dot{\gamma}(t)}{|\dot{\gamma}(t)|}$ and $\ell(\gamma)$ is the length of $\gamma$ (see \Cref{Def: PlaneCurve}). 
%Since $\gamma$ is closed, i.e., for the interval $I=[0,1]$, we have that $\gamma(0)=\gamma(1)$,
Moreover, for any $x=\gamma(t)$ we denote the tangent vector $\tau(x):=\tau(t)$ whenever this exists, and its $\frac{\pi}{2}$-counterclockwise rotation by $\eta(x)$. A clear geometric presentation of the above is with:

\begin{equation}\label{eq: Main1a}
 G_\delta(\gamma)= \ell(\gamma) + \delta\int_{\gamma(\RZ)}\int_{\gamma(\RZ)}\frac{|\eta(x)-\eta(x')|^p}{|x-x'|^{1+sp}}d\mathcal{H}^1(x)d\mathcal{H}^1(x'),
 \end{equation}
 \end{subequations}
\noindent Equation ~\eqref{eq: Main1a} suggests a nonlocal oscillation of the normal vector hence describing a notion of geometric curvature for a curve.

The framework of \eqref{Main} models the study of soap films with the features of thickness and volume enclosure. The dimensional setting presents a most natural notion of span without recourse to the Harrison-Pugh model. Furthermore, the analysis does not require much sophisticated machinery or esoteric sets of geometric measure theory. Moreover, a key feature of embeddability for finite energy curves is shown to correspond to a \textit{non-collapsing} phenomenon.

Our first key result is on well-posedness of the minimisation problem \eqref{Main}:

\begin{theorem}[Existence of minimizers]\label{thm:Existence} There exists a solution $\gamma\in \mathcal{E}$ to the minimization problem (\ref{Main}).
\end{theorem}

A key ingredient in the proof is a regularity result for curves with finite energy. Furthermore, we obtain an embedding result in a suitable fractional Sobolev space which corresponds to a non-collapsing feature in our model. More precisely, we show

\begin{theorem}[Non-collapsing for finite-energy curves]\label{Thm: CharFiniteEnergy}
Let $\gamma\in C^{0,1}(\RZ, \R^2)$ be a constant-speed curve and $s\in (0,1),\, p\in (1,\infty)$ with $sp > 1$ and $\delta>0$. Then $G_\delta <\infty$ implies that $\gamma$ is embedded and $\gamma\in W^{1+s,p}(\RZ,\R^2)$.
\end{theorem}

\Cref{Thm: CharFiniteEnergy} can be viewed as a geometric Sobolev embedding, a technique with roots in knot theory which characterises Sobolev spaces using geometric functionals. A slightly more general converse, \Cref{Prop: CharFiniteEnergy}, is established. Furthermore, we give a counter-example (see \Cref{Counterexample}) to \Cref{Thm: CharFiniteEnergy}  in the case of parametric region $sp<1$. In the case of $sp=1$, We make the interesting observation of connection to the so-called Mobius energy in knot theory via the Ishizeki-Nagasawa decomposition (see e.g., \cite{IN14Mob}).

Our next main result is a partial $\Gamma$-convergence of $G_\delta$ as $\delta\rightarrow 0$. This gives the convergence of our model to the capillarity model considered in Maggi \textit{et al.} \cites{MSS19,KMS20}.  The result is stated in the subclass $\bar{\mathcal{E}}_R$ (see definition in \eqref{Def:Class E-R}) of curves in $\mathcal{E}$ obtained from \textit{rearrangements} (see \Cref{Def: Rearr}) which admit smooth embedded approximations. To state the result, we first define $L_\eps$ on $W^{1,\infty}$ as
\[
L_\eps(\gamma):= \begin{cases}
\|\dot{\gamma}\|_\infty,  & \gamma \in \bar{\mathcal{E}}_R,\\
 +\infty,& \text{otherwise}.
\end{cases}
\]

\begin{theorem}[Partial $\Gamma-$convergence as $\delta\rightarrow 0$]\label{Gconv}
Let $\gamma\in W^{1,\infty}$. Then
\begin{enumerate}
\item[(i)] for any sequence $\gamma_\delta$ converging to $\gamma$ uniformly in $W^{1,\infty}$ we have that
\[
L_\eps(\gamma)\le \liminf\limits_{\delta}G_\delta(\gamma_\delta)
\]
\item[(ii)] if $\gamma\in \bar{\mathcal{E}}_R$ there exists a sequence $(\gamma_\delta)$ converging to $\gamma$ uniformly in $W^{1,\infty}$ and satisfying
\[
\limsup\limits_{\delta}G_\delta(\gamma_\delta)\le L_\eps(\gamma).
\]
\end{enumerate}
%\[
%\Gamma-\lim_{\delta\rightarrow 0}G_\delta (\gamma)= L(\gamma)
%\]
%with respect to weak star convergence in $W^{1,\infty}$.
\end{theorem}

The proof of \Cref{Gconv} uses a technical lemma (\Cref{Prop: Rearr}) on rearrangement of Lipschitz curves with winding number of $0$ or $1$. We note that such curves can in general have self-intersections. The rearrangement decreases length and admits smooth embedded approximations. The approximating sequence is reconstructed via some perturbation and area fixing operations to guarantee inclusion in the class $\mathcal{E}$. A direct consequence of the Theorem is the convergence of minimizers of $G_\delta$:

\begin{corollary}[Convergence of minimizers]\label{Cor:ConvMinimizers}
Let $\{\gamma_\delta\}$ be a sequence of minimizers of $G_\delta$ such that
\[
\gamma_\delta \rightarrow \gamma \quad \text{ uniformly in } W^{1,\infty} \text{ as } \delta\rightarrow 0.
\]
Then $\gamma$ is a minimizer of $L_\eps$.
\end{corollary}

We call minimizers $\gamma$ of $L_\eps$  \textit{capillary curves} for natural reason of positive area enclosure. We remark that by restricting to the subclass $\mathcal{E}_R$, the minimization is unchanged due to an equivalence result on the minimization in the classes $\mathcal{E}_R$ and $\mathcal{E}$ as being the same (see \Cref{Prop: ClassEquiv}).

The next main result is on $\Gamma$-convergence of $L_\eps$ as $\eps\rightarrow 0$. We consider the class $\bar{\mathcal{E}}_0$ (see \Cref{Def: Class-0}) of closed lipschitz curves $\gamma$ satisfying $\mathrm{E}1$ and $\omega(\gamma,x)=0$ for all $x\notin\gamma(\RZ)$ and define the function $L_0$ on $W^{1,\infty}$ as
\begin{equation}\label{Eqn:Final limit}
L_0(\gamma):=\begin{cases}
\|\dot{\gamma}\|_\infty,\quad \gamma\in \bar{\mathcal{E}}_0,\\
+\infty,\quad \text{otherwise}.
\end{cases}
\end{equation}

\begin{theorem}[$\Gamma-$convergence to Plateau problem]\label{Thm: eps Gconv}
$L_\eps$ $\Gamma$-converges to $L_0$ as $\eps\rightarrow 0$ with respect to uniform convergence. That is, for any $\gamma\in W^{1,\infty}$ we have that
\begin{enumerate}
\item[(i)] for any sequence $\gamma_\eps$ converging to $\gamma$ uniformly in $W^{1,\infty}$ as $\eps\rightarrow 0$ we have that
\[
L_0(\gamma)\le \liminf\limits_{\eps\rightarrow 0}L_{\eps}(\gamma_{\eps})
\]
\item[(ii)] there exists a sequence $(\gamma_\eps)$ converging to $\gamma$ uniformly in $W^{1,\infty}$ as $\eps\rightarrow 0$ and satisfying
\[
\limsup\limits_{\eps\rightarrow 0}L_\eps(\gamma_\eps)\le L_0(\gamma).
\]
\end{enumerate}
\end{theorem}
In particular, \Cref{Thm: eps Gconv} gives convergence to minimizers of $L_0$, we refer to these as \textit{Plateau curves}. This seals up the sum on \Cref{Prob:Plateau}. To see this, we prove:

\begin{proposition}\label{Thm:Concluding}
Given $\gamma\in \bar{\mathcal{E}}_0$, let $K:=\gamma(\RZ)$. Then
\[
L_0(\gamma) \ge  2\H^1(K).
\]
On the other hand, for every compact connected set $K\subset \R^2$ there exists $\gamma_K\in \bar{\mathcal{E}}_0$ such that $K=\gamma(\RZ)$ and
\[
L_0(\gamma_K)= 2\H^1(K).
\]
\end{proposition}
And finally say,

\begin{theorem}\label{Cor:Concluding}
Given a set $K_0$ which is a minimizer for \Cref{Prob:Plateau}, there exists a rectifiable curve $\gamma_0\in \bar{\mathcal{E}}_0$ such that 
\begin{equation}\label{Eq:norm-length Equiv.}
L_0(\gamma_0) =  2\H^1(K_0),
\end{equation}
and $\gamma_0$ is a minimizer for the functional $L_0$. On the other hand, given a minimizer $\gamma\in \bar{\mathcal{E}}_0$ of $L_0$, it's image $K:=\gamma(\RZ)$ is a solution to \Cref{Prob:Plateau}.
\end{theorem}

\subsection{Review of the Maggi \textit{et al.} model}
The study in Maggi \textit{et al.} \cite{KMS20} introduces a new perspective on modeling soap films at equilibrium, treating them not merely as surfaces but as regions of small total volume. This framework is within the context of capillarity problems with homotopic spanning conditions, which allows for a more nuanced understanding of the physical behavior of soap films. They prove, among other things, the existence of \textit{generalized minimizers} for the capillarity problems (see \cite[Theorem~1.4]{KMS20}) which converge to a solution of the classical Plateau's problem (in the spirit of Harrison-Pugh) recovered in the vanishing volume limit (see \cite[Theorem~1.9]{KMS20}). 

We remark that by resorting to an enclosed volume rather than thickness in the discourse of \cite{KMS20}, the minimizers of the Capillary problem though possessing regions of positive volume, due to the spanning condition, in general fail to have uniformly positive thickness.
The effect of this is the occurrence of collapsed regions in the equilibrium configuration of minimizers -- the \textit{collapsing phenomenon}. A collapsed section marks a region of discontinued observation of soap films. Further analysis of these regions become of essential interest to construe the model (see \cites{KMS20smoothness, KMS20collapsing}). The build-up of this exciting approach of Maggi \textit{et al.} uses  exotic sets from geometric measure theory and delicate machinery from the calculus of variations via the direct method approach, a groundwork inspired by Del~Lellis \textit{et al.} \cites{DlGM17,DlDrG19} and De~Philipis \textit{et al.} \cite{DpDrG16}.

A possible idea in addressing the issue of collapsed regions can be seen from a heuristic insight into the collapsing phenomenon in the Maggi \textit{et al.} model, which is this: Soap films can admit singularities (see  \cites{almgren1976geometry,taylor1976structure}), and when these are present, with the introduction of small but positive volume, minimizers collapse in order to achieve stability and optimize energy. This can be linked to the fact that soap films possess thickness  which shows in what is sometimes observed in physical soap film experiments (see \cite{boys1959soap}, \cites{almgren1976geometry,taylor1976structure}). The motivation then is to consider a simple model which factors the essential feature of thickness of soap films whilst still retaining the capillarity formulation of the problem. This essentially is the inspiration of the present work.

\subsection{Structure of paper}
%We now give an outline of the paper. 
In \Cref{Sec:Prelim} we collect useful notions and results related to closed curves in the plane, and of fractional Sobolev spaces. The proof of \Cref{thm:Existence} is given in \Cref{Sec:Existence}. This uses a regularization property, \Cref{Lem:FinEnergy}, of finite energy curves. The section ends with proof of \Cref{Thm: CharFiniteEnergy} and a converse result \Cref{Prop: CharFiniteEnergy}. We also give an analysis of the parameter regions $sp<1$ and $sp=1$ which entails a counter-example in the case of the former, and connection to geometric knot energies for the latter.

\Cref{Sec:PartialGamma} contains the partial $\Gamma$-convergence result \Cref{Gconv}, and \Cref{Cor:ConvMinimizers} on the convergence of minimizers. These give a sense of our model converging to the capillary model of Maggi \textit{et al.} \cites{MSS19,KMS20}. The proofs use a rearrangement result \Cref{Prop: Rearr} and a class fixing procedure \Cref{Prop: RearrSeq}.

Finally, in \Cref{Sec:Gamma} we prove the $\Gamma$-convergence result \Cref{Thm: eps Gconv} which gives convergence of the Capillarity model to the Plateau problem. Moreover, we obtain a solution to \Cref{Prob:Plateau}. This becomes clear from  proofs of \Cref{Thm:Concluding} and \Cref{Cor:Concluding} presented in the same section as well.

\section{Preliminaries}\label{Sec:Prelim}
%We are interested in a more general class of closed curves that are not necessarily simple. The following notions are essential for the definition of such. 
We denote $C^{0,1}(\RZ)$ as the space of closed Lipschitz curves in the plane.

\begin{definition}[Plane curve]\label{Def: PlaneCurve}
Let $\gamma \in C^{0,1}(\RZ)$. Then $\gamma$ is regular if $|\dot{\gamma}(t)|\ne 0$ for almost every $t\in \RZ$. It is said to have constant speed if there exists $C>0$ such that $|\dot{\gamma}(s)|=C$ for almost every $s\in \RZ$. We define length of $\gamma$,
\[
\ell(\gamma)=\sup _{\left(t_{i}\right)_{i}} \sum_{i=1}^{k}\left|\gamma\left(t_{i}\right)-\gamma\left(t_{i-1}\right)\right|=\int\limits_{\RZ}|\dot{\gamma}|dt
\]
where the supremum is taken over all partitions $0=t_{0}<t_{i}<\cdots<t_{k}=1$.\\
A curve is rectifiable if it has finite length.
\end{definition}
\begin{remark}\label{Rmk:Orientation}
We note that if $\gamma$ is a rectifiable curve then for $\H^1$-a.e. $x\in\gamma(\RZ)$ there exists an approximate tangent vector $\tau(x)$ to $\gamma(\RZ)$ at $x$, $\gamma^{-1}(\{x\})$ is finite, and $\dot{\gamma}(t)$ exists and is nonzero for every $t\in \gamma^{-1}(\{x\})$. Furthermore, $\gamma(\RZ)$ can be oriented by a collection of all such tangent vectors. (see, e.g., \cite[(3.2.19)]{Fed69} or \cite{Falconer1986})
\end{remark}
\begin{remark}\label{Rmk:speed-length}
We note that If $\gamma\in C^{0,1}(\RZ)$ is a curve with constant speed, the Lipschitz constant of $\gamma$ coincides with the length $\ell(\gamma)$. In this case, the functional $G$ can be re-written in the following form
\[
G_\delta(\gamma)= \ell(\gamma) + \delta\ell(\gamma)^{2-p}\iint\limits_{\RZ\times \RZ}\frac{\left|\dot{\gamma}(t)-\dot{\gamma}(t')\right|^p}{|\gamma(t)-\gamma(t')|^{1+sp}}dtdt'.
\] 
\end{remark}

The following is a well known reparametrization result for curves which have finite length (see, for instance, \cite[Remark 3.2]{AlbertiOttolini17}).

\begin{lemma}\label{Lem:ContReparam.}
If $\gamma:\RZ\rightarrow \R^2$ is continuous and $\ell:=\ell(\gamma)<\infty$. Then there exists $\sigma:[0,1]\rightarrow \RZ$ an increasing function such that $\gamma\circ\sigma$ is Lipschitz and has constant speed $\ell$. 
\end{lemma}
\begin{definition}\label{Def: WindingNum}
Given $\gamma : \RZ\rightarrow \mathbb{R}^2\equiv \mathbb{C}$ a closed  Lipschitz plane curve. For any $x\in \mathbb{R}^2\setminus \gamma(I)$, we define the winding number
\begin{equation}\label{eq: deg}
\omega(\gamma, x)=\frac{1}{2\pi i}\int_{\gamma(\RZ)}\frac{d\zeta}{\zeta-x} =\frac{1}{2\pi i}\int_{\RZ}\frac{\dot{\gamma}(t)}{\gamma(t)-x}dt. 
\end{equation}
\end{definition}
The following proposition recalls the additivity property of the winding number under concatenation of curves.

\begin{proposition}\cite[Theorem 7.2]{StewartTall83}
Let $\gamma_1$ and $\gamma_2$ be continuous closed curves in $\mathbb{C}\setminus \{0\}$ such that the end point of $\gamma_1$ is the start of $\gamma_2$. Then
\[
\omega(\gamma_1\star\gamma_2,0)=\omega(\gamma_1,0)+\omega(\gamma_2,0).
\]
\end{proposition}

\begin{remark}\label{Rmk:CurveAssumption} For the rest of this paper we only consider curves $\gamma\in C^{0,1}(\RZ)$ with $\omega(\gamma,x)\in \{0,1\}$ for all $x\ne \gamma(\RZ)$. In this case we define the interior set of $\gamma$
\[
\mathtt{int}(\gamma) =\{x\in \mathbb{R}^2 : \omega(\gamma,x)=1\},
\]
\noindent We shall denote by $|\mathtt{int}(\gamma)|$ the Lebesgue measure of $\mathtt{int}(\gamma)$. A simple calculation gives that
\[
|\mathtt{int}(\gamma)|=\int\limits_{\R^2\setminus \gamma(\RZ)}\omega(\gamma, x)dx=\frac{1}{2\pi i}\int\limits_{\R^2\setminus \gamma(\RZ)}\int\limits_{\RZ}\frac{\dot{\gamma}(t)}{\gamma(t)-x}dtdx
\]
Moreover, for $\gamma(t)=\left(\gamma_1(t),\gamma_2(t)\right)$, by Green's theorem we have
\begin{eqnarray*}
|\mathtt{int}(\gamma)|&=&-\int_{\RZ}\gamma_2(t)\cdot\dot{\gamma}_1(t)dt = \int_{\RZ}\gamma_1(t)\cdot\dot{\gamma}_2(t)dt\\
 &=& \frac{1}{2}\int_{\RZ}\gamma_1(t)\dot{\gamma}_2(t)-\gamma_2(t)\dot{\gamma}_1(t)dt = \int_{\RZ}\frac{\gamma(t)^\perp\cdot\dot{\gamma}(t)}{2}.
\end{eqnarray*}
\end{remark}

\begin{proposition}\cite[Corollary~2]{JB91}\label{Prop: DegCont}
Let $\gamma_j:\RZ\rightarrow \R^2$ for all $j\in \mathbb{N}$ and $\gamma:\RZ\rightarrow \R^2$ be curves such that $\gamma_j \rightarrow \gamma$ uniformly.
%and $\dot{\gamma}_j\rightharpoonup \dot{\gamma}$ weakly in $L^1$
 Then for any $x\in \R^2\setminus\gamma(\RZ)$ and for all $j$ large enough,
$\omega(\gamma_j,x)= \omega(\gamma,x)$.
\end{proposition}

\noindent The following is a standard result of topology (see e.g., \cite[page 67]{BerensteinGay} or \cite[page 321]{BergGost88}).

\begin{proposition}\label{Lem: IntNum}
Let $\gamma\in C^1(\RZ)$ be a closed  regular curve and let $\phi:[0,1]\rightarrow \R^2$ be a $C^1$ curve with $\dot{\phi}\ne 0$ such that $\phi$ is transversal to $\gamma$ with endpoints $\phi(0)=x_0$ and $\phi(1)=x_1$ both not in $\gamma(\RZ)$, then
\begin{equation}\label{Eq:IntNum}
\omega(\gamma,x_1)-\omega(\gamma, x_0)= \text{\textit{``intersection number} of $\gamma$ and $\phi$''}
\end{equation}
where the intersection number of $\gamma$ and $\phi$ is the number of intersections of $\gamma$ and $\phi$ counted with orientations, that is, $+1$ for all $(t,s)$ such that $\gamma(t)=\phi(s)$ and $\{(\dot{\gamma}(t),\dot{\phi}(s))\}$ define a positive orientation and $-1$ if they have a negative orientation.
\end{proposition}

\subsection{Fractional Sobolev spaces}
We briefly recall the fractional Sobolev spaces. We refer to reference \cite{DnPV12} for a handy exposition on the topic.
\begin{definition}
Let $\alpha=(\alpha_1,\ldots, \alpha_n)$ having natural components $\alpha_i$ as a standard multi-index notation with $|\alpha|=\alpha_1+\cdots+\alpha_n$. We define the partial derivative
\[
D^\alpha u=\frac{\partial^{|\alpha|}u}{\partial x_1^{\alpha_1}\cdots\partial x_n^{\alpha_n}}
\]
Let $p\in [1,\infty),\, k,n\in \mathbb{N},\, s\in (0,1)$ and $\Omega\subset\R^n$ be open. We define the fractional Sobolev space $W^{k+s,p}(\Omega)$ by 
\[
W^{k+s,p}(\Omega)=\left\{u\!\in\! W^{k,p}(\Omega): \int_\Omega\int_\Omega\frac{|D^\alpha u(x)-D^\alpha u(y)|^p}{|x-y|^{n+sp}}dxdy < \infty \,\forall \alpha \text{ s.t. } |\alpha|=k\right\}
\]
We denote by 
\[
\|D^\alpha u\|_{W^{s,p}}^p:=\int_\Omega\int_\Omega\frac{|D^\alpha u(x)-D^\alpha u(y)|^p}{|x-y|^{n+sp}}dxdy
\]
Endowed with the norm
\[
\|u\|^p_{W^{k+s,p}}:=\|u\|^p_{W^{k,p}}+\sum\limits_{|\alpha|=k}^{}\|D^\alpha u\|^p_{W^{s,p}},
\] 
$W^{k+s,p}(\Omega)$ is a Banach space. 
\end{definition}
The following is an analogue of the embedding theorem of Sobolev spaces (See e.g., \cite[Chapter 4]{AdamsFournier03} and \cite[Theorem 4.54]{demengel2012functional})

\begin{theorem}[Fractional Sobolev embedding] Let $\Omega \subset \mathbb{R}^{n}$ be a Lipschitz domain, $n, p \in \mathbb{N}, 1 \leq p<\infty$ and $s \in(0,1)$: 
\begin{enumerate}
\item[(a)] If $sp>n$, then there is a constant $C=C(n, p, s, \Omega)$ such that for any $u \in W^{s,p}(\Omega)$ we have
\[
|| u||_{C^{0, \alpha}(\Omega)} \leq C\|u\|_{W^{s,p}}
\]
with $\alpha=\frac{(sp-n)}{p}$\\
Furthermore, we have for $k \in \mathbb{N}$
\[
\|u\|_{C^{k, \alpha}(\Omega)} \leq C\|u\|_{W^{k+s, p}}
\]
We always consider the continuous representative in the Sobolev classes.
\item[(b)] If $sp<n$, then there is a constant $C=C(n, p, s, \Omega)$ such that for any $u \in W^{s,p}(\Omega)$ we have
\[
\|u\|_{{L^q}(\Omega)} \leq C \|u\|_{W^{s,p}}, \text{ for every } q\in [p,p^*]
\]
where $p^*=\frac{np}{n-sp}$ is the fractional critical exponent. Moreover, we have the compact embedding
\[
W^{s,p}(\Omega)\hookrightarrow L^q(\Omega),\text{ for every } q\in [1,p^*).
\]
\end{enumerate}
\end{theorem}
\section{Existence and Non-collapsing of Minimizers}\label{Sec:Existence}

\noindent In this section, we present the proofs of \Cref{thm:Existence} and \Cref{Thm: CharFiniteEnergy}. We denote, for every $s\in (0,1)$, $1\le p\le \infty$ and $\gamma\in C^{0,1}(\RZ)$,
\begin{equation}
G(\gamma):= \int_{\RZ}\int_{\RZ}\frac{\left|\tau(t)-\tau(t')\right|^p}{|\gamma(t)-\gamma(t')|^{1+sp}}|\dot{\gamma}(t)||\dot{\gamma}(t')|dtdt'. 
 \end{equation}

Clearly, $G$ is the nonlocal term in $G_\delta$ and contributes to analytic and geometric properties of the functional $G_\delta$.  

\subsection{Proof of \Cref{thm:Existence}: Existence of minimizers}

We first prove the following elementary lemmas for the proof of \Cref{thm:Existence}.

\begin{lemma}\label{Lem:FinEnergy}
Let $\gamma\in C^{0,1}(\RZ)$ be a constant-speed curve. If $G\!<\infty$ then $\gamma\in W^{1+s,p}(\RZ,\R^2)$.
\end{lemma}
\begin{proof}
Using that $\gamma$ has constant speed, we have that (see \Cref{Rmk:speed-length})
\begin{eqnarray*}
\infty>\; G &= & \ell(\gamma)^{2-p}\iint\limits_{\RZ\times \RZ}\frac{\left|\dot{\gamma}(t)-\dot{\gamma}(t')\right|^p}{|\gamma(t)-\gamma(t')|^{1+sp}}dtdt'\\
&\ge & \ell(\gamma)^{1-(1+s)p}\iint\limits_{\RZ\times \RZ}\frac{\left|\dot{\gamma}(t)-\dot{\gamma}(t')\right|^p}{|t-t'|^{1+sp}}dtdt.
\end{eqnarray*}
It follows that $\gamma\in W^{1+s,p}(\RZ, \R^2)$.
\end{proof}

\begin{lemma}\label{Lem: continuity}
Let $\gamma_j:\RZ\rightarrow \R^2$ and $\gamma:\RZ\rightarrow \R^2$ be constant-speed curves for all $j\in \mathbb{N}$. Suppose $\gamma_j \rightarrow \gamma$ in $W^{1+s,p}$. Then $\ell(\gamma_j) \rightarrow \ell(\gamma):=|\dot{\gamma}|$.
\end{lemma}
\begin{proof}
By the Sobolev embedding of $W^{s,p}$ in $L^1$ we obtain, using the triangle inequality, that

\[
|\ell(\gamma_j)-\ell(\gamma)|\le \int_{\RZ}\Big||\dot{\gamma}_j|-|\dot{\gamma}|\Big|\le \int_{\RZ}|\dot{\gamma}_j-\dot{\gamma}|\rightarrow 0
\]
as $j\rightarrow \infty$.
\end{proof}

\begin{lemma}\label{Prop: AreaCont}
Let $\gamma_j:\RZ\rightarrow \R^2$ and $\gamma:\RZ\rightarrow \R^2$ be constant-speed curves for all $j\in \mathbb{N}$. Suppose $\gamma_j \rightarrow \gamma$ uniformly. Then 
$\lim\limits_{j\rightarrow\infty}|\mathtt{int}(\gamma_j)|= |\mathtt{int}(\gamma)|$.
\end{lemma}
\begin{proof}
By uniform convergence, the sequence $\gamma_j$ is uniformly bounded and hence the enclosures $\mathtt{int}(\gamma_j)$ are contained in some open disk. By the dominated convergence theorem we have that
\[
|\mathtt{int}(\gamma_j)|=\int \mathbf{1}_{\mathtt{int}(\gamma_j)}\rightarrow\int \mathbf{1}_{\mathtt{int}(\gamma)}= |\mathtt{int}(\gamma)|
\]
as $j\rightarrow \infty$
\end{proof}

\begin{lemma}\label{Lem: lsc}
Let $\gamma_j:\RZ\rightarrow \R^2$ and $\gamma:\RZ\rightarrow \R^2$ be constant-speed curves for all $j\in \mathbb{N}$. Suppose $\gamma_j \rightarrow \gamma$ uniformly and $\dot{\gamma}_j\rightarrow \dot{\gamma}$ pointwise almost everywhere. Then 
\[
\liminf G_\delta(\gamma_j)\ge G_\delta(\gamma).
\]
\end{lemma}
\begin{proof}
This is a direct application of Fatou's lemma. %where we have used the fact that the product of weakly convergent and uniformly convergent sequences gives a pointwise convergent sequence 
\end{proof}

\begin{proof}[Proof of \Cref{thm:Existence}]
We proceed in the following steps:
 \begin{enumerate}
\item[(i)] $\mathcal{E}$ is nonempty: % since for any unit-speed circle $\beta$ of area $\eps$ one can deform using area-fixing methods in such a way that $\beta \supset \{x_1,\ldots, x_N\}$ and then reparametrize to arclength. 
We think this claim is natural and so omit the proof. A simple heuristic idea is to consider the circle of area $\eps$ deformed using area-preserving Lipschitz operations to contain the points $\{x_1,\ldots, x_N\}$. %See figure [...]

\item[(ii)] Uniform bounds:
 Let $\left\lbrace \gamma_j \right\rbrace_j \subset \mathcal{E}$ be a minimizing  sequence with
 \[
G_\delta(\gamma_j)\le C<+\infty \quad \text{for all } j\in \mathbb{N}, \]
for some $C>0$ (we can take $2\ell(\gamma_0)$). 
Clearly $\ell_j:=\ell(\gamma_j)\le C$ for all $j\in \mathbb{N}$. Moreover $\{\dot{\gamma_j}\}$ is uniformly bounded in $W^{s,p}(\RZ)$, this follows from estimate in the proof of \Cref{Lem:FinEnergy}. Furthermore, there exists $C'>0$ such $\ell_j\ge C'$ for all $j\in \mathbb{N}$. Indeed, since $\gamma_j$ contains at least two distinct points (e.g., $x_1,x_2$ as in condition $\mathrm{E}1$),
\[
\ell_j\ge \ell\left(\gamma_j\big|_{\{x_1, x_2\}
}\right)\ge \|x_1 - x_2\| \quad \text{for all } j\in \mathbb{N},
\]
where $\gamma_j\big|_{\{x_1, x_2\}}$ is the arc of $\gamma_j$ joining the points $x_1, x_2$.
\item[(iii)] Convergences:
Now, since $\gamma_j$ fixes $N$ points we have that the sequence $\gamma_j$ and their enclosures are contained in some open disk. This contains every $\mathtt{int}(\gamma_j)$. By the parametrization we have uniformly Lipschitz property and hence by the Ascoli-Arzela theorem, we have that $\gamma_j \rightarrow \gamma$ uniformly (up to a subsequence) for some Lipschitz curve $\gamma$. Furthermore we have that $\dot{\gamma}_j\rightarrow \dot{\gamma}$ weakly-star in $L^\infty$. Also by the uniform bound on $\dot{\gamma}_j $ in $W^{s,p}(\RZ)$ and  the compact embedding of $W^{s,p}(\RZ)$ in $L^1$
% for some $q\in [1,p^*)$, where $p^*=\frac{p}{1-sp}$ is the  critical fractional Sobolev exponent of $p$,
  we have that $\dot{\gamma_j}\rightarrow \dot{\gamma}$ in $L^1$ and hence pointwise almost everywhere (up to a subsequence). Using the continuity of the length functional as in \Cref{Lem: continuity}, we have that $\ell_j \rightarrow \ell:=\ell(\gamma)$ and hence  $|\dot{\gamma}|=$ constant (since $|\dot{\gamma}_j|$ is a sequence of constant functions).
\item[(iv)] Closure: We show $\gamma \in \mathcal{E}$. 
Indeed, from the continuity of $\omega(\cdot,x)$ with respect to uniform convergence as in \Cref{Prop: DegCont} we obtain that $\omega(\gamma,x)\in \{0,1\}$ for any $x\notin \gamma(\RZ)$. By (iii) above we readily have that $\gamma$ has constant speed. Also we establish by \Cref{Prop: AreaCont} that $|\mathtt{int}(\gamma)|=\eps$. Finally, the uniform  convergence  gives that $\gamma \supset \{x_1,\ldots, x_N\}$.

%%%%%---------------%%%%%%%%%%%%%%%%%%%%%%%%%%ALTERNATIVELY, even if $||f_n-f||_p\rightarrow 0$ does not imply that $f_n\rightarrow f(x)$ a.e., it is possible to find a subsequence $\{f_{n_k}\}$ such that $f_{n_k}(x)\rightarrow f(x)$ a.e. [SEE Harold E. Krogstad ``some informal notes about $L^p$-spaces and convergence"]. So we can have the result up to a subsequence. KENNEDY LET'S MOVE ON.

%%%%%%-----------%%%%%%%%%%%%%%%%%%%%%%%%%%%%%
 \end{enumerate}
We conclude the proof by the lower semicontinuity result \Cref{Lem: lsc}.\qedhere
 
\end{proof}

 \begin{remark}
In general, minimizers can exhibit pathological behaviours such as self-intersections. However, the characterization of curves with finite-energy, gives that minimizers have improved regularity in the $sp>1$ case. 
 \end{remark}

\subsection{Finite-energy curves and proof of \Cref{Thm: CharFiniteEnergy}}

In what follows, we show that in the parameter region corresponding to $sp>1$, we recover nice analytic and geometric features of curves from the finiteness of the energy. In particular, we obtain that a rectifiable curve $\gamma$ with finite $G$ (hence and necessarily for finite $G_\delta$) is embedded.\\[5pt]

\noindent
We first prove the following crucial bi-Lipschitz estimate for injective $C^1$ curves followed by a simple lemma which captures the nature of tangent vectors at self-intersection points.

\begin{lemma}\label{lem: BiLip}
Let $\gamma\in C^1(\RZ)$ be an injective and regular curve. Then $\gamma$ is an embedding. Moreover, $\gamma$ is bi-Lipschitz, i.e., there exists $M>0$ such that 
\begin{equation}\label{eq:biLip}
|\gamma(s)-\gamma(t)|\ge M|s-t|, \quad \text{ for all } s,t \in \RZ.
\end{equation}
\end{lemma}

\begin{proof}
A direct application of the inverse function theorem gives that $\gamma$ is an embedding and $\gamma(\RZ)$ is a submanifold of $\R^2$. Now define
$F:\RZ\times \RZ \rightarrow (0,\infty)$ such that
\[
F(s,t)=\begin{cases}
\frac{|\gamma(s)-\gamma(t)|}{|s-t|}, & s\ne t,\\
|\dot{\gamma}(t)|, & s=t.
\end{cases}
\]
Since $\gamma \in C^1(\RZ,\R^2)$, $F$ is a continuous. Also, since $\gamma$ is injective and regular, we have that $F$ is strictly positive. Furthermore, since the domain is compact the minimum is attained and strictly positive. Denoting by $M:=\min\limits_{\RZ\times \RZ}F$, 
the inequality (\ref{eq:biLip}) follows.\end{proof}

%\begin{lemma}\label{Lem: Collinear}
%Let $\gamma: \RZ \rightarrow \R^2$ be a $C^1$ constant-speed curve with $\omega(\gamma,x)\in \{0,1\}$ for all $x\notin \gamma(\RZ)$. Then every multiple point  has velocity vectors that are collinear
%\end{lemma}

\begin{lemma}\label{Lem: Collinear}
Let $\gamma: \RZ \rightarrow \R^2$ be a constant-speed curve with $\omega(\gamma,x)\in \{0,1\}$ for all $x\notin \gamma(\RZ)$. Let $x\in \gamma(\RZ)$ with $\gamma^{-1}{(\{x\})}=\{t_1,t_2\}$ and suppose  $\dot{\gamma}(t_1)$ and $\dot{\gamma}(t_2)$ are well-defined. Then the tangent vectors  $\dot{\gamma}(t_1)$ and $\dot{\gamma}(t_2)$ are parallel and with opposite directions i.e.  $\dot{\gamma}(t_1) =- \dot{\gamma}(t_2)$.
\end{lemma}
\begin{proof}The proof is essentially \cite[Lemma 4.3]{DPP20} if we take $u=I$.
%Suppose $\dot{\gamma}(t_1)$ and $\dot{\gamma}(t_2)$ are not collinear then we can find some neighbourhood $V_x$ of $x$ and a path $\beta$ transversal to $\gamma$ and intersecting it two times (see Fig. -). Now by Lemma \ref{Lem: IntNum}, we obtain a difference of two in the winding numbers of the endpoints of $\beta$. This is a contradiction to the assumption on the winding numbers.
\end{proof}

%\begin{lemma}\label{SelfTanVectors}
%Let $\gamma: \RZ \rightarrow \R^2$ be a $C^1$ constant-speed curve with $\omega(\gamma,x)\in \{0,1\}$ for all $x\notin \gamma(\RZ)$. Then if $\gamma$ is not injective, there exist $t_1,\, t_2\in \RZ$, $t_1\ne t_2$ such that $\gamma(t_1)=\gamma(t_2)$ and  $\dot{\gamma}(t_1)\ne\dot{\gamma}(t_2)$. In particular, $\dot{\gamma}(t_1)=-\dot{\gamma}(t_2)$.
%\end{lemma}

\begin{lemma}\label{SelfTanVectors}
Let $\gamma: \RZ \rightarrow \R^2$ be a constant-speed curve with $\omega(\gamma,x)\in \{0,1\}$ for all $x\notin \gamma(\RZ)$. Let $x\in \gamma(\RZ)$ such that $\gamma^{-1}{(\{x\})}=\{t_1,\ldots,t_m\}\subset \RZ$ with $m\ge 2$, and suppose  $\dot{\gamma}(t_i)$ is well-defined for all $1\le i\le m$. Then there exist $1\le k\ne l\le m$, such that  $\dot{\gamma}(t_k)\ne\dot{\gamma}(t_l)$.
\end{lemma}

\begin{proof}
Without loss of generality, suppose $\gamma$ has unit-speed parametrization. We claim there exist $1\le k\ne l\le m$ such that $\dot{\gamma}(t_k)\ne \dot{\gamma}(t_l)$. Suppose otherwise, that is, for every $1\le i,j\le m$, $\dot{\gamma}(t_i)=\dot{\gamma}(t_j)$. Let $\phi :\RZ\rightarrow \R^2$ be a $C^1$ curve transversal to $\gamma$ through $x$ and having endpoints $\phi(0)=x_0$ and $\phi(1)=x_1$. By \Cref{Lem: IntNum} we immediately obtain a winding number difference of at least two between the endpoints. This is certainly not possible under the assumption on the winding number.
\end{proof}

\noindent We now present:
%\begin{theorem}[Non-collapsing effect of finite-energy curves]\label{thm: CharFiniteEnergy}
%Let $\gamma\in C^{0,1}(\RZ, \R^2)$ be a constant-speed curve and $s\in (0,1),\, p\in (1,\infty)$ with $sp > 1$. Then $G<\infty$ implies that $\gamma$ is embedded and $\gamma\in W^{1+s,p}(\RZ,\R^2)$.
%\end{theorem}

\begin{proof}[Proof of \Cref{Thm: CharFiniteEnergy}]
Let $\gamma \in C^{0,1}(\RZ)$ such that $G_\delta (\gamma) < \infty$. By \Cref{Lem:FinEnergy} we have that $\gamma\in W^{1+s,p}(\RZ, \R^2)$. We now show $\gamma$ is injective in $\RZ$.  We proceed by a contradiction. Suppose $\gamma(t_1)=\gamma(t_2)$ for some  $t_1\ne t_2$. By \Cref{SelfTanVectors} can suppose that $\dot{\gamma}(t_1)\ne\dot{\gamma}(t_2)$. Since $sp>1$, by the Sobolev embedding of $W^{s,p}(\RZ)$ into $C^{0,\alpha}(\RZ)$ for some $\alpha$, we have that $\gamma\in C^1$. Consider the Taylor's expansion:
\begin{eqnarray*}
\gamma(t')=&\gamma(t_1)+\dot{\gamma}(t_1)(t'-t_1)+o(|t'-t_1|^{1+\alpha}),\\
\gamma(t)=&\gamma(t_2)+\dot{\gamma}(t_2)(t-t_2)+o(|t-t_2|^{1+\alpha}),
\end{eqnarray*}
and obtain that
\[
|\gamma(t')-\gamma(t)|\le C\left(|t'-t_1|+|t-t_2|\right)
\]
for some $C>0$. By the continuity of $\dot{\gamma}$ we consider $0<\rho<1$ such that
\begin{eqnarray*}
|\dot{\gamma}(t)-\dot{\gamma}(t')|\ge c>0,\qquad\quad
\end{eqnarray*}   
for some constant $c$, in the neighbourhood, $I_\rho:=\{|t-t_1|<\rho\}\times\{|t'-t_2|<\rho\}$ of $t_1,\, t_2$.

Now we have,
 \begin{eqnarray*}
G(\gamma)&\!=\!&\ell(\gamma)^{2-p}\iint\limits_{\RZ\times \RZ}\frac{\left|\dot{\gamma}(t)-\dot{\gamma}(t')\right|^p}{|\gamma(t)-\gamma(t')|^{1+sp}}|dtdt'\\
&\!\ge\!& \frac{c^p}{C^{1+sp}}\ell(\gamma)^{2-p}\iint_{I_\rho}\frac{1}{\left(|t'-t_1|+|t-t_2|\right)^{1+sp}}dtdt'\\
&\!=\!& \frac{c^p}{C^{1+sp}}\ell(\gamma)^{2-p}\int_{\{|t|<\rho\}}\!\!\int_{\{|t'|<\rho\}}\frac{1}{\left(|t'|+|t|\right)^{1+sp}}dtdt'.
\end{eqnarray*}
Passing to polar coordinates one easily checks that this diverges which gives a contradiction. Hence $\gamma$ is injective and we conclude using \Cref{lem: BiLip}.
\end{proof}

The following is a converse to \Cref{Thm: CharFiniteEnergy} in a slightly more general form.
\begin{proposition}[\textit{Converse} to \Cref{Thm: CharFiniteEnergy}]\label{Prop: CharFiniteEnergy}
Let $s\in (0,1),\, p\in (1,\infty)$ with $sp > 1$ and $\delta>0$. Let $\gamma\in W^{1+s,p}(\RZ,\R^2)$ be injective and having constant speed. Then $G_\delta <\infty$.
\end{proposition}

\begin{proof}
Using 
%the elementary inequality,
%\[
%\left|\frac{a}{|a|}-\frac{b}{|b|}\right|\le 2\frac{|a-b|}{|a|}, \quad \text{for any } a,b\in \R
%\]
%and 
the bi-Lipschitz estimate \Cref{lem: BiLip}, we obtain that
\begin{eqnarray*}
G_\delta (\gamma) &=& \ell(\gamma) + \delta \ell(\gamma)^{2-p}\int\int_{\RZ^2}\frac{|\dot{\gamma}(t)-\dot{\gamma}(t')|^p}{|\gamma(t)-\gamma(t')|^{1+sp}}dtdt'\\
&\le & \ell(\gamma) + \delta\frac{ \ell(\gamma)^{2-p}}{M^{1+sp}}\int\int_{\RZ^2}\frac{|\dot{\gamma}(t)-\dot{\gamma}(t')|^p}{|t-t'|^{1+sp}}dtdt' \quad < \infty,
\end{eqnarray*}
where $M$ is the bi-Lipschitz constant in  (\ref{eq:biLip}).
\end{proof}

\begin{remark}
\Cref{Thm: CharFiniteEnergy} presents a sort of geometric Sobolev embedding which is in the spirit of characterizing Sobolev spaces using geometric functionals (see e.g., \cite{Kol15Geom}, \cite{BoNg06} and \cite{BrNg20}). This allows for analytic information from geometric functionals which encode data involving steric interactions. Furthermore, combined with the \cref{lem: BiLip}, we see that the functional $G_\delta$ suggests a notion of knottedness or a measure of \textit{entangledness} of finite-energy curves.
\end{remark}

\begin{remark}[The $sp<1$ case]
The embedding theorem for finite-energy curves suggests a natural domain to consider in the space of curves: parameter region corresponding to $sp>1$. In fact, we see that in the case $sp<1$, we do not only lose embeddability but also the non-collapsing property as the following example illustrates.

\begin{example}[Counter-example in the $sp<1$ case]\label{Counterexample}
We consider the following example of a self-intersecting curve with finite-energy in the neighbourhood of an intersection.

Let 
\[
\begin{cases}
\gamma(y)=(y,-1), & y\in(-\frac{1}{10},\frac{1}{10}),\\
\tilde{\gamma}(x)=(\sin x,\cos x), & x\in (\pi-\alpha,\pi+\alpha), \; \alpha>0.
\end{cases}
\]
We have that, $\gamma(0)=\tilde{\gamma}(\pi)$ and $\dot{\gamma}(0)=-\tilde{\gamma}'(\pi)$ and obtain

\[
\int\limits_{-\frac{1}{10}}^{\frac{1}{10}}\int\limits_{\pi-\alpha}^{\pi+\alpha}\frac{|(1-\cos x, \sin x)|^{p}}{|(y-\sin x, -1-\cos x)|^{1+sp}} = \int\limits_{-\frac{1}{10}}^{\frac{1}{10}}\int\limits_{\pi-\alpha}^{\pi+\alpha}\frac{((1-\cos x)^2+\sin^2 x)^{\frac{p}{2}}}{((y-\sin x)^2+ (-1-\cos x)^2)^{\frac{1+sp}{2}}}
\]
%\\[-25pt]
\begin{eqnarray*}
&\le&\int\limits_{-\frac{1}{10}}^{\frac{1}{10}}\int\limits_{\pi-\alpha}^{\pi+\alpha}\frac{5^{\frac{p}{2}}}{\left((y+(x-\pi)+O((x-\pi)^3))^2+(\frac{1}{2}(x-\pi)^2+O((x-\pi)^3)^2)\right)^{\frac{1+sp}{2}}}\\
&\le&\int\limits_{-\frac{1}{10}}^{\frac{1}{10}}\int\limits_{\pi-\alpha}^{\pi+\alpha}\frac{5^{\frac{p}{2}}}{\left(y^2+2y(x-\pi)+(x-\pi)^2+(\frac{1}{4}(x-\pi)^4)\right)^{\frac{1+sp}{2}}}\\
&\le&\int\limits_{-\frac{1}{10}}^{\frac{1}{10}}\int\limits_{\pi-\alpha}^{\pi+\alpha}\frac{5^{\frac{p}{2}}}{\left(y^2+2y(x-\pi)+(x-\pi)^2)\right)^{\frac{1+sp}{2}}}=\int\limits_{-\frac{1}{10}}^{\frac{1}{10}}\int\limits_{-\alpha}^{\alpha}\frac{5^{\frac{p}{2}}}{\left(y^2+2xy+x^2)\right)^{\frac{1+sp}{2}}}<\infty
\end{eqnarray*}

\end{example}

Hence, we see that for the $sp<1$ case, we expect collapsing effect (as in the example above). However, we remark that due to the regularity property, this is not the case for minimizers. 
%Indeed, from the Euler-Lagrange equation we obtain regularity property of minimizers from which we deduce non-collapsing.
\end{remark}

\begin{remark}[The $sp=1$ case]
We make the important observation that the case of the parameter region corresponding to  $sp=1$ with $p=2$, our nonlocal term $G$ is actually the first integrand in the Mobius energy via the Ishizeki-Nagasawa decomposition (see \cite{IN14Mob})
\[
E_{\mathrm{Mob}}(\gamma)=E_{1}(\gamma)+E_{2}(\gamma)+4
\]
where
\[
E_{1}(\gamma):=\iint_{(\mathbb{R} / \mathbb{Z})^{2}} \frac{|\tau(x)-\tau(y)|^{2}}{2|\gamma(x)-\gamma(y)|^{2}}\left|\gamma^{\prime}(x) \| \gamma^{\prime}(y)\right| \mathrm{d} x \mathrm{d} y,
\]
$\tau=\frac{\gamma^{\prime}}{\left|\gamma^{\prime}\right|}$ and

\begin{eqnarray*}
E_{2}(\gamma)\!\!&:=&\!\!\!\! \iint_{(\mathbb{R}  \mathbb{Z})^{2}} \frac{2}{|\gamma(x)-\gamma(y)|^{2}}  \operatorname{det}\left(\!\!\begin{array}{cc}\langle\tau(x), \tau(y)\rangle & \langle(\gamma(x)-\gamma(y)), \tau(x)\rangle \\ \langle(\gamma(x)-\gamma(y)), \tau(y)\rangle & |\gamma(x)-\gamma(y)|^{2}\end{array}\!\!\right)\\[5pt]
&& \hspace{8cm}\times\left|\gamma^{\prime}(x) \| \gamma^{\prime}(y)\right| \mathrm{d} x \mathrm{d} y.
\end{eqnarray*}
Various studies such as obtaining explicit variational formulae, proving existence and uniqueness of minimizers and regularity of critical points of functionals with a Mobius energy have been carried out (see e.g., \cite{HeZ00},\cite{FHe94}, \cite{BR13OHara}), in particular, the decomposition above (see \cite{IN15Mob}).
\end{remark}

%\section{$\Gamma$-Convergence result as $\delta\rightarrow 0$}
\section{Convergence to the capillarity model}\label{Sec:PartialGamma}

The present section contains a convergence result which shows that letting $\delta$ tend to zero, we obtain a limit curve $\gamma$ of the sequence of minimisers $\gamma_\delta$ of the functionals $G_\delta$.
More precisely, we prove a partial $\Gamma$-convergence result for the sequence of functionals $G_\delta$ in the (closure of  the) class $\mathcal{E}$. This requires the construction of an appropriate subclass $\mathcal{E}_R$ of \textit{rearrangements} which admit good approximations. We refer to the limit curve $\gamma$ as a capillary curve. 
We also write the closure  $\bar{\mathcal{E}}$ {(with respect to weak star convergence in $W^{1,\infty}$)} of $\mathcal{E}$ which gives the
%closure of $\mathcal{E}$ with respect to weak star convergence in $W^{1,\infty}$.
class of all Lipschitz curves $\gamma$ that satisfy $\mathrm{E}1 - \mathrm{E}3$ but not necessarily $\mathrm{E}4$.
\vspace{5pt}

We introduce the following terminology.
\begin{definition}\label{Def: Rearr}
Let $\gamma\in C^{0,1}(\RZ)$. We say $\tilde{\gamma}\in C^{0,1}$ is a rearrangement of $\gamma$ if there exists $\sigma:\RZ\rightarrow \RZ$ which is Borel and injective such that $\tilde{\gamma}=\gamma\circ \sigma$
\end{definition}

\noindent
We now begin with the following technical lemma due to G. Alberti \cite{Alb2021} which helps in defining an appropriate class for the  partial $\Gamma$-convergence result.

\begin{lemma}\cite{Alb2021}\label{Prop: Rearr}
Let $\gamma\in W^{1,\infty}(\RZ,\R^2)$ be a closed curve with $\omega(\gamma,x)\in \{0,1\}$ for all $x\notin \gamma(\RZ)$. Then there exists a rearrangement $\tilde{\gamma}$ of $\gamma$ satisfying
\begin{enumerate}
\item[R1.] $\tilde{\gamma}(\RZ)=\gamma(\RZ)$ and $\#(\tilde{\gamma}^{-1}(x))\le \#(\gamma^{-1}(x))$ for all $x\in \gamma(\RZ)$. In particular, $\ell(\tilde{\gamma})\le \ell(\gamma)$;
\item[R2.] $\omega(\tilde{\gamma},x)=\omega(\gamma,x)\in \{0,1\}$ for every $x\notin \gamma(\RZ)$;
\item[R3.] there exists a sequence $(\gamma_n)$ of smooth embedded regular closed curves oriented counter-clockwise such that $\gamma_n\rightarrow \tilde{\gamma}$ uniformly;
\item[R4.] $\|\dot{\gamma}_n\|_\infty\rightarrow \|\dot{\tilde{\gamma}}\|_\infty$.
\end{enumerate}
\end{lemma}

%\begin{remark}
%The proof is essentially showing a reparametrization which gives self-tangency in place of transversal self-crossings and a perturbation which transforms the self-tangent curve into a simple curve. Then we conclude with a standard smoothing (piecewise) approximation.
%\end{remark}

\noindent We show in the sequel that for every $\gamma\in \mathcal{E}$, the corresponding rearrangement approximating sequence  in \Cref{Prop: Rearr} satisfies (up to some modification) the admissibility conditions $(\mathrm{E}1)-(\mathrm{E}4)$.
\begin{proposition}\label{Prop: RearrSeq}
Let $\gamma\in \mathcal{E}$. Then the approximating sequence $(\gamma_n)$ corresponding to the rearrangement $\tilde{\gamma}$ of $\gamma$ as  in \Cref{Prop: Rearr}, up to some modifications, satisfies items $(\mathrm{E}1)-(\mathrm{E}4)$.
\end{proposition}

\begin{proof}
Clearly, from item $R3$, $\omega(\gamma_n,x)\in \{0,1\}$ for every $x\notin \gamma_n(\RZ)$ for all $n$. Since $\gamma_n$ is regular, up to a reparametrization, we can suppose the curves $\gamma_n$ have constant speed.\\
\textbf{Step 1} (Uniform area): By the continuity of the area functional with respect to the uniform convergence of the sequence (see \Cref{Prop: AreaCont}), we have that
\[
\eps_n:=|\mathtt{int}(\gamma_n)|\longrightarrow \eps:=|\mathtt{int}(\gamma)|
\]
Rescaling by a factor $\lambda_n=\sqrt{\frac{\eps}{\eps_n}}$, we can just suppose $|\mathtt{int}(\gamma_n)|=\eps$ for all $n\ge 1$. (Note that $\lambda_n\rightarrow 1$ as $n\rightarrow \infty$).\\[5pt]
\textbf{Step 2} (Perturbation): Now we perturb our curves in a slight way to accommodate the fixed points $\{x_1,\ldots,x_N\}$. More precisely, we modify the curves $\gamma_n$ slightly by a small perturbation near the fixed points and subsequently apply some area-fixing operations. 
%The idea is to use local variations, essentially along the lines of \cite[Lemma 29.13]{Mag12}. 
We proceed as follows:\\[3pt]
Step 2(a): For every $1\le i\le N$, let $\zeta_i>0$ be small enough (and to be chosen). 
%let $\delta_i=\delta(\zeta_i)>0$. 
%and $m\in \N$ such that $B_\delta(x_i)\cap \gamma_n(\RZ)\ne \emptyset$ and $B_\delta(x_j)\cap B_\delta(x_k)=\emptyset$ for all $n\ge m, 1\le i,j,k\le N, \; j\ne k$. 
For any $n$ large enough, consider point $x_i'\in \gamma_n(\RZ)$ which minimizes $d(x_i,\gamma_n(\RZ))$.
Denote $m_i=\min d(x_i,x_i')$. If $m_i=0$ then we are done. Suppose $m_i>0$, since $\gamma_n$ is smooth, we find $\eta_i=\eta_i(\zeta_i)$ such that denoting by $\partial(B_{\eta_i}(x_i')\cap \gamma_n(\RZ)):=\{x_{i,a}',x_{i,b}'\}$ the two boundary points closest to $x_i$, we obtain the smooth {triangular spike} $T_i=T(x_{i,a}',x_{i,b}',x_i)$ with base given by arc in $B_{\eta_i}(x_i')\cap \gamma_n(\RZ)$ joining $\{x_{i,a}',x_{i,b}'\}$, vertices $\{x_{i,a}',x_{i,b}',x_i\}$ and area $\zeta_i$. The perturbed curve (not relabelled) consists of the curve $T_i$ replacing the sub-curve $B_{\eta_i}(x_i')\cap \gamma_n(\RZ)$ in $\gamma_n$ (see Figure \ref{Fig:Perturbation}).\\[3pt]
\noindent
Step 2(b): The final stage in the perturbation is fixing the area. This uses the fact that $\zeta_i$ can be chosen as small as possible. Using that $\gamma_n$ is injective and $|\mathtt{int}(\gamma_n)|=\eps>0$, we find a suitable open disk $B_i$ not containing any of the points $x_1,\ldots,x_N$ such that $\H^1(\gamma_n(\RZ)\cap B_i)>0$. Denoting by $\sigma_i:=\sigma_0\left(\mathtt{int}(\gamma_n),B_i\right)$ as in Lemma~17.21 in \cite{Mag12}, we choose $\zeta_i\in(0,\sigma_i)$ and obtain an area-fixing variation which ``moves" $\gamma_n(\RZ)\cap B_i$ in order to increase or decrease area of $\mathtt{int}(\gamma_n)\cap B_i$ by $\zeta_i$. More precisely, if $x_i\in \mathtt{int}(\gamma_n)$ for any $1\le i\le N$, then the area loss in the construction in step $(2a)$ is recovered by an area increase from the local variation. Conversely, if $x_i\notin \mathtt{int}(\gamma_n)$ for any $1\le i\le N$, the area gained is corrected by an area decrease from the variation. Moreover, the local variation of \cite{Mag12}*{Lemma~17.21} ensures the properties of $\gamma_n$ are preserved throughout the operation.
\begin{figure}[!hbtp]
\centering
\includegraphics[width=0.5\linewidth]{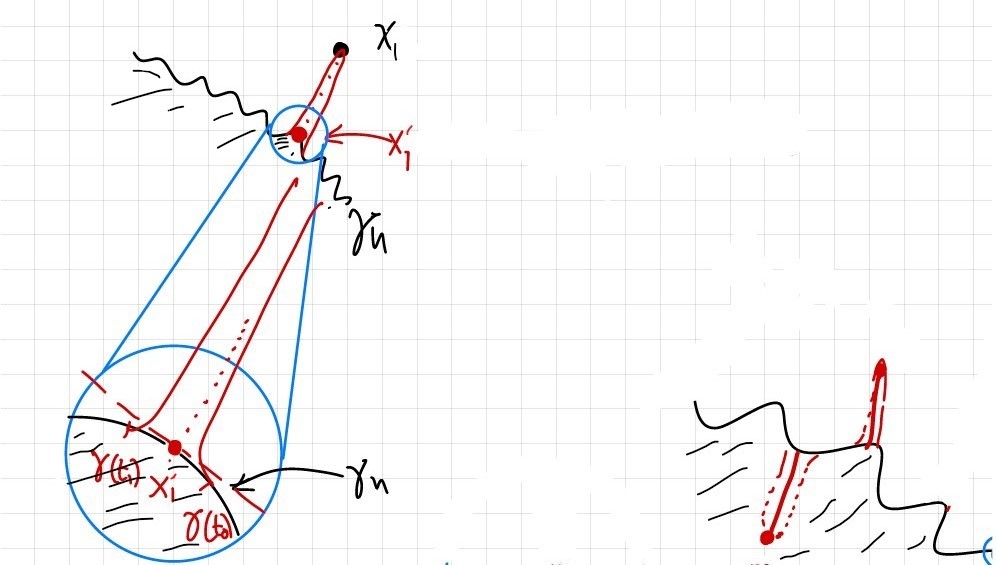}
\caption{Perturbation of curve near fixed points}
\label{Fig:Perturbation}
\end{figure}
\noindent\\
\textbf{Step 3}: Finally, up to constant speed reparametrization, we obtain that the modified sequence of curves $\gamma_n$ (without relabelling) satisfies the conditions $(\mathrm{E}1)-(\mathrm{E}4)$.
\end{proof}
\vspace{15pt}
For any $\delta>0$, we extend $G_\delta$ to the whole $W^{1,\infty}(\RZ)$ by $+\infty$, that is,
\[
G_\delta(\gamma)= \begin{cases}
 \|\dot{\gamma}\|_\infty + \delta \ell(\gamma)^{2-p}\int\limits_{\RZ}\int\limits_{\RZ}\frac{|\dot{\gamma}(t)-\dot{\gamma}(t')|^p}{|\gamma(t)-\gamma(t')|^{1+sp}}dtdt',& \gamma \in \mathcal{E},\\
 +\infty,& \text{otherwise}.
\end{cases}
\]

We define the class
\begin{equation}\label{Def:Class E-R}
\bar{\mathcal{E}}_R :=\left\lbrace\tilde{\gamma
} \text{ obtained from applying \Cref{Prop: Rearr}  to }\gamma\in \bar{\mathcal{E}}\right\rbrace
\end{equation}

\begin{remark}
We observe that the class $\bar{\mathcal{E}}_R$ is strictly contained in $\bar{\mathcal{E}}$. Indeed, consider the quadrifolium curve, see \Cref{Fig: Quadrifolium}.

\tikzset{               %new code
    redarrows/.style={postaction={decorate},decoration={markings,mark=at position 0.1 with {\arrow{>}}},decoration={markings,mark=at position 0.4 with {\arrow{>}}},decoration={markings,mark=at position -0.4 with {\arrow{>}}},decoration={markings,mark=at position -0.1 with {\arrow{>}}}},
    redarrowsh/.style={postaction={decorate},decoration={markings,mark=at position 0.1 with {\arrow{>}}},decoration={markings,mark=at position 0.4 with {\arrow{>}}},decoration={markings,mark=at position 1 with {\arrow{>}}}}}
\begin{figure}[!htb]
\centering
%\begin{subfigure}[t]{0.4\textwidth}
%\centering
\begin{tikzpicture}
\draw[domain=0:360, scale=1,samples=500,redarrows] plot (\x: {1.5*sin(2*\x)});
\draw[color=blue] (1.4,1) node {$1$};
\draw[color=blue] (-1.4,1) node {$2$};
\draw[color=blue] (-1.4,-1) node {$3$};
\draw[color=blue] (1.4,-1) node {$4$};
\end{tikzpicture}
\caption{The Quadrifolium traversed in the order $1$-$4$-$3$-$2$ is an example of a curve in $\bar{\mathcal{E}}$ which is not in $\bar{\mathcal{E}}_R$.
%A rearrangement for this is traversing in the order $1$-$2$-$3$-$4$.
}
%\end{subfigure} %
%~
%\begin{subfigure}[t]{0.4\textwidth}
%\centering
%\begin{tikzpicture}
%\draw[domain=0:360, scale=1,samples=500,redarrowsh] plot(\x: {1.5*cos(3*\x)});
%\end{tikzpicture}
%\caption{Trefolium}
%\end{subfigure}
%\caption{Examples of curves in $\bar{\mathcal{E}}$ that are not in $\bar{\mathcal{E}}_R$}
\label{Fig: Quadrifolium}
\end{figure}
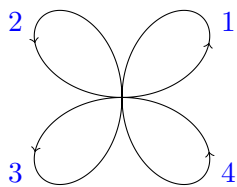
\end{remark}
We now proceed first with the following compactness lemma which gives the domain of the $\Gamma-$limit as $W^{1,\infty}\cap \bar{\mathcal{E}}$.

\begin{lemma}[Compactness]\label{Lem: Compactness}
Let $(\delta_n)$ be a sequence with $\delta_n\rightarrow 0$ as $n\rightarrow\infty$ and let $(\gamma_n)\subset \mathcal{E}$ also be given such that $G_{\delta_n}(\gamma_n)<M$ for some $M<\infty$. Then there exist a subsequence $\gamma_{n_k}$ and $\gamma\in W^{1,\infty}\cap\bar{\mathcal{E}}$ such that $\gamma_{n_k}\rightarrow \gamma$ uniformly in $W^{1,\infty}$.
\end{lemma}
%\begin{remark}
%Moreover, $\gamma$ is possibly self-intersecting but without self-crossings.
%\end{remark}
\begin{proof}
To simplify notations we represent $\delta_n,\, \gamma_n$ and $G_{\delta_n}$ by $\delta, \gamma_\delta$ and $G_\delta$ respectively. Using the uniform bound condition $G_\delta(\gamma_\delta)<M$ we obtain by  \Cref{Lem:FinEnergy} that $\gamma_\delta\in W^{1+s,p}$. Furthermore, by the formulation of $G_\delta$, we obtain a uniform bound on the derivatives, precisely,
\begin{equation}\label{eq:equibd}
\|\dot{\gamma}_\delta\|_{\infty} \le M \quad \text{ for all } \delta\in (0,1).
\end{equation}
Moreover by (\ref{eq:equibd}) we have that the sequence is equi-Lipschitz. Using the condition $\mathrm{E}3$, let $x_0$ be any of the fixed points and for any $\delta\in (0,1)$, take $t_\delta$ such that $\gamma_\delta(t_\delta)=x_0$. Then by the Lipschitz property and (\ref{eq:equibd}), we have  that
\[
|\gamma_\delta(t)-x_0|\le M|t-t_\delta|\le M,\qquad \text{ for all } t\in \RZ
\]
and hence obtain that
\[
\sup_\delta|\gamma_\delta|\le M + |x_0|.
\] 
Therefore, by the Ascoli-Arzela theorem, we obtain that $\gamma_\delta\rightarrow \gamma$ (up to a subsequence) uniformly for some $\gamma\in W^{1,\infty}$. Now, by uniform convergence we have that $\gamma(\RZ)\supset \{x_1,\ldots,x_N\}$ and, by \Cref{Prop: DegCont}, that $\omega(\gamma,x)\in \{0,1\}$ for every $x\notin \gamma(\RZ)$. Also, by \Cref{Prop: AreaCont}, we have that $|\mathtt{int}(\gamma)|=\eps$. Hence the conclusion follows that $\gamma\in W^{1,\infty}\cap \bar{\mathcal{E}}$.
%Finally, we have that $\gamma$ is possibly self-intersecting but does not admit self-crossings as a uniform limit of simple curves (since self-crossing is invariant under the convergence).
\end{proof}

\begin{remark}
We now present the proof of \Cref{Gconv} but make the observation that the $\limsup$ inequality is established essentially in the subclass $\bar{\mathcal{E}}_R$ of curves admitting nice approximations. However, we will later see that this is sufficient to obtain the convergence of minimizers.
\end{remark}

\begin{proof}[Proof of \Cref{Gconv}]
The lower bound inequality: Suppose $\gamma\in W^{1,\infty}$ and $\gamma_\delta\rightarrow \gamma$ uniformly in $W^{1,\infty}$. We can suppose $(\gamma_\delta)\subset \mathcal{E}$ otherwise the lower bound inequality readily holds. Then clearly $\gamma\in \bar{\mathcal{E}}$ and by the lower semicontinuity of $L_\eps$ and the inequality $L_\eps\le G_\delta$ (recall that $L_\eps(\gamma_\delta)=\ell(\gamma_\delta)$ because $\gamma_\delta$ has constant speed) for every $\delta>0$, the inequality
\[
L_\eps(\gamma)\le \liminf\limits_{\delta}G_\delta(\gamma_\delta)
\]
follows.\\ 
To prove the upper bound (UB) inequality we consider $\gamma\in W^{1,\infty}$. The UB inequality holds trivially for $\gamma\notin \bar{\mathcal{E}}_R$. Suppose $\gamma\in \bar{\mathcal{E}}_R\subset\bar{\mathcal{E}}$. By definition we have that $\omega(\gamma
,x)\in \{0,1\}$ for every $x\notin \gamma(\RZ)$ and $\gamma(\RZ)\supset\{x_1,\ldots,x_N\}$, and also we have that $|\mathtt{int}(\gamma)|=\eps$. Now, we consider the sequence $(\gamma_n)$ of smooth regular closed curves in \Cref{Prop: Rearr} corresponding to the rearrangement $\tilde{\gamma}$ of $\gamma$. By \Cref{Prop: RearrSeq} we can suppose that $\gamma_n$ satisfies conditions $(\mathrm{E}1)-(\mathrm{E}4)$ for every $n$, hence in
$\mathcal{E}$.
Now using that $\|\dot{\gamma}_n\|_\infty\rightarrow \|\dot{\tilde{\gamma}}\|_\infty$ which means that $L_\eps(\gamma_n)\rightarrow L_\eps(\tilde{\gamma})$ and by a suitable choice of subsequence $n=n(\delta)$ depending on $\delta$ for all $\delta\in(0,1)$ such that $\delta G(\gamma_n)\rightarrow 0$ as $\delta\rightarrow 0$,
%(e.g. appropriating constant sequence terms for decreasing values of $\delta$), 
we obtain that
\begin{equation}\label{Ineq: limsup}
\limsup\limits_{\delta}G_\delta(\gamma_\delta)\le L_\eps(\gamma).
\end{equation}
and this concludes the proof.
\end{proof}

\noindent As a  consequence of the partial  $\Gamma$-convergence result, we establish the result on convergence of minimizers.
\begin{comment}
\begin{corollary}[Convergence of minimizers]
Let $\{\gamma_\delta\}$ be a sequence of minimizers of $G_\delta$ such that
\[
\gamma_\delta \rightarrow \gamma \quad \text{ uniformly in } W^{1,\infty}.
\]
Then $\gamma$ is a minimizer of $L$.
\end{corollary}
\end{comment}

\begin{proof}[Proof of \Cref{Cor:ConvMinimizers}]
By the $\Gamma$-liminf inequality we have that
\[
L_\eps(\gamma)\le \liminf\limits_{\delta\rightarrow 0} G_\delta(\gamma_\delta)=\liminf\limits_{\delta\rightarrow 0}\inf\limits_{W^{1,\infty}}G_\delta.
\]
Let $\beta\in W^{1,\infty}$ be arbitrary. We want to prove that $L_\eps(\beta)\ge L_\eps(\gamma)$. If $\beta\notin \bar{\mathcal{E}}_R$, then we are done. Hence we suppose $\beta\in \bar{\mathcal{E}}_R$ and let $\tilde{\beta}\in \bar{\mathcal{E}}_R$ its corresponding rearrangement. Let $\beta_\delta$ be a recovery sequence of $\tilde{\beta}$, then
\[
\liminf\limits_{\delta\rightarrow 0}\inf\limits_{W^{1,\infty}}G_\delta\le\limsup\limits_{\delta\rightarrow 0}G_\delta(\beta_\delta)\le L_\eps(\tilde{\beta})\le L_\eps(\beta)
\]
and we are done.
\end{proof}

\begin{definition}\label{Def: classes}
We define $\tilde{\mathcal{E}}$ as the class of all curves $\gamma$ such that
\begin{enumerate}
\item[(i)] $\gamma$ is continuous and has finite length,
\item[(ii)] $\gamma(\RZ)\supset \{x_1,\ldots, x_N\}$,
\item[(iii)] $\omega(\gamma,x)\in \{0,1\}$ for every $x\notin \gamma(\RZ)$,
\item[(iv)] $|\mathsf{int}(\gamma)|=\eps$.
\end{enumerate}
Observe that curves $\gamma$ in $\tilde{\mathcal{E}}$ need not have constant speed. It is not difficult to see that $\bar{\mathcal{E}}_R\subset \tilde{\mathcal{E}}$.
\end{definition}

This section is concluded with an equivalence result which reconciles the $\Gamma$-limit functional and the minimisation problem with the standard setting of minimizing the length functional in the more general class $\tilde{\mathcal{E}}$. This essentially says the partial $\Gamma$-convergence result, although given in the class $\bar{\mathcal{E}}_R\subset \tilde{\mathcal{E}}$, suffices for our variational investigation.

\begin{proposition}\label{Prop: ClassEquiv}
Let $\bar{\mathcal{E}}$, $\tilde{\mathcal{E}}$ and  $\bar{\mathcal{E}}_R$ be as in \Cref{Def: classes} and (\ref{Def:Class E-R}). Then
\[
\inf\limits_{\bar{\mathcal{E}}_R}L_\eps(\gamma)=\inf\limits_{\bar{\mathcal{E}}}L_\eps(\gamma)=\inf\limits_{\tilde{\mathcal{E}}}\ell(\gamma)
\]
\end{proposition}

\begin{proof}
Recall the inclusions $\tilde{\mathcal{E}}\supset\bar{\mathcal{E}}\supsetneq\bar{\mathcal{E}}_R$. We begin with the first equality. Clearly $\inf\limits_{\bar{\mathcal{E}}}L_\eps(\gamma)\le \inf\limits_{\bar{\mathcal{E}}_R}L_\eps(\gamma)$. To show the converse inequality, we proceed thus:
Let $\gamma'\in \bar{\mathcal{E}}$ be arbitrary. Using \Cref{Prop: Rearr}, we let $\tilde{\gamma}'\in \bar{\mathcal{E}}_R$ be a corresponding rearrangement with $\ell(\tilde{\gamma}')\le \ell(\gamma')$. Without loss of generality we can suppose $\tilde{\gamma}'$ has constant speed. We obtain, recalling \Cref{Rmk:speed-length}, that
\[
\inf\limits_{\bar{\mathcal{E}}_R}L_\eps(\gamma)\le L_\eps(\tilde{\gamma}')=\ell(\tilde{\gamma}')\le\ell(\gamma')\le L_\eps(\gamma').
\]
Since $\gamma'$ is arbitrary the desired converse inequality follows.

For the second equality: since $\bar{\mathcal{E}}\subset\tilde{\mathcal{E}}$, we have that 
\[
\inf\limits_{\tilde{\mathcal{E}}}\ell(\gamma)\le  \inf\limits_{\bar{\mathcal{E}}}\ell(\gamma)\le \inf\limits_{\bar{\mathcal{E}}}L_\eps(\gamma).
\]
Now, for any $\gamma\in \tilde{\mathcal{E}}$, by \Cref{Lem:ContReparam.} there exists a reparametrization $\tilde{\gamma}$ which is Lipschitz, has constant speed and satisfies 
\[
\ell(\gamma)=\ell(\tilde{\gamma})=L_\eps(\tilde{\gamma})\ge \inf\limits_{\bar{\mathcal{E}}}L_\eps(\gamma).
\]
Hence it follows that 
\[
\inf\limits_{\tilde{\mathcal{E}}}\ell(\gamma)\ge\inf\limits_{\bar{\mathcal{E}}}L_\eps(\gamma).
\]
This concludes the proof.
\end{proof}

%\section{$\Gamma$-Convergence result as $\eps\rightarrow 0$}
\section{Plateau problem as capillarity limit}\label{Sec:Gamma}

In this concluding section we prove a $\Gamma$-convergence result which is used to obtain a solution to \Cref{Prob:Plateau}.
We first explicitly define the classes:
\begin{definition}\label{Def: Class-0}
For any $\eps>0$, we denote the class $\bar{\mathcal{E}}_\eps$ as $\bar{\mathcal{E}}$ from \Cref{Def: classes} and further 
introduce the class $\bar{\mathcal{E}}_0$ of all curves $\gamma$ satisfying
\begin{enumerate}
\item[(i)] $\gamma\in W^{1,\infty}(\RZ)$ closed, Lipschitz curve,
\item[(ii)] $\gamma(\RZ)\supset \{x_1,\ldots,x_N\}$,
\item[(iii)] $\omega(\gamma,x)=0$ for all $x\notin\gamma(\RZ)$. In particular, $|\mathsf{int}(\gamma)|=0$.
\end{enumerate}
\end{definition}
We proceed to recall the function $L_0$ on $W^{1,\infty}$:
\[
L_0(\gamma):=\begin{cases}
\|\dot{\gamma}\|_\infty,\quad \gamma\in \bar{\mathcal{E}}_0\\
+\infty,\quad \text{otherwise}.
\end{cases}
\]
We now prove the $\Gamma$-convergence result \Cref{Thm: eps Gconv} following the compactness result below.

\begin{lemma}[Compactness]\label{Lem: eps-Compactness}
Let $(\eps_n)$ be a sequence with $\eps_n\rightarrow 0$ as $n\rightarrow\infty$ and let $(\gamma_n)\subset \mathcal{E}(\eps_n)$ also be given such that $L_{\eps_n}(\gamma_n)<M$ for some $M<\infty$. Then there exist a subsequence $\gamma_{n_k}$ and $\gamma\in W^{1,\infty}\cap\bar{\mathcal{E}}_0$ such that $\gamma_{n_k}\rightarrow \gamma$ uniformly in $W^{1,\infty}$.
\end{lemma}

\begin{proof}
The uniform convergence result (up to a subsequence, of course) follows by direct arguments as in the proof of \Cref{Lem: Compactness}. We claim $\gamma\in \bar{\mathcal{E}}_0$. Clearly, by uniform convergence we have that $\gamma(\RZ)\supset \{x_1,\ldots,x_N\}$. Furthermore, 
\[
\eps_n:=|\mathtt{int}(\gamma_n)|\rightarrow |\mathtt{int}(\gamma)|=0
\]
which gives that $\omega(\gamma,x)=0$ for all $x\notin \gamma(\RZ)$. Hence the claim follows.
\end{proof}

\begin{proof}[Proof of \Cref{Thm: eps Gconv}]
We prove the liminf inequality. Let $\gamma_\eps\rightarrow \gamma$ uniformly in $W^{1,\infty}$ as $\eps\rightarrow 0$. We can suppose $\liminf\limits_{\eps\rightarrow 0}L_\eps(\gamma_\eps)<+\infty$, otherwise there is nothing to prove. Now take a minimizing subsequence $(\eps_n)$ such that
\[
\lim\limits_{\eps_n\rightarrow 0}L_{\eps_n}(\gamma_{\eps_n})=\liminf\limits_{\eps\rightarrow 0}L_{\eps}(\gamma_{\eps}).
\]
Then $\gamma_{\eps_n}\in \bar{\mathcal{E}}_{\eps_n}$ and by the lower semicontinuity of the $L^\infty$-norm we obtain that
\[
L_0(\gamma) \le \liminf\limits_{\eps_n\rightarrow 0}L_{\eps_n}(\gamma_{\eps_n})
\]
Clearly, from the compactness result \Cref{Lem: eps-Compactness}, we have that $\gamma\in \bar{\mathcal{E}}_0$ and the liminf inequality follows.

For the limsup inequality, we construct a sequence $(\gamma_\eps)\subset \bar{\mathcal{E}}$ with $\gamma_\eps\rightarrow \gamma$ uniformly in $W^{1,\infty}$ such that 
\[\limsup\limits_{\eps}L_\eps(\gamma_\eps)\le L_0(\gamma).\]
Clearly, if $\gamma\notin \bar{\mathcal{E}}_0$, the inequality holds trivially. So suppose $\gamma\in \bar{\mathcal{E}}_0$. For some $\eps>0$ (small enough), we denote $\lambda:=\eps^{\frac{1}{4}}$. Let $\tilde{\gamma}:[0,1]\rightarrow \R^2$ be an arclength parametrized circle with radius $\sqrt{\frac{\eps}{\pi}}$. We define 
\[
\gamma_\eps(t):=\begin{cases}
\tilde{\gamma}(\frac{t}{\lambda}),& t\in[0,\lambda],\\
\gamma(\frac{t-\lambda}{1-\lambda}),& t\in[\lambda,1].
\end{cases}
\]
Clearly, using the additivity property of the winding number, it is not difficult to see that $\gamma_\eps\in \bar{\mathcal{E}}_\eps$. We note that $\|\dot{\tilde{\gamma}}\|_\infty=2\sqrt{\pi\eps}$ and thus obtain that 

\begin{eqnarray*}
\|\dot{\gamma}_\eps\|_\infty&\le& \max\left\lbrace\frac{1}{1-\lambda}\|\dot{\gamma}\|_\infty,\frac{1}{\lambda}\|\dot{\tilde{\gamma}}\|_\infty\right\rbrace\\
&\le&\max\left\lbrace\frac{1}{1-\eps^{\frac{1}{4}}}\|\dot{\gamma}\|_\infty,2\sqrt{\pi}\eps^{\frac{1}{4}}\right\rbrace,
\end{eqnarray*}
hence 
\[
\limsup\limits_{\eps}\|\dot{\gamma}_\eps\|_\infty \le \|\dot{\gamma}\|_\infty=L_0(\gamma)
\]
and the limsup inequality follows.
\end{proof}

In the sequel we present proofs to \Cref{Thm:Concluding} and \Cref{Cor:Concluding} which reconcile the results of the present section with solutions to \Cref{Prob:Plateau}.

\begin{proof}[Proof of \Cref{Thm:Concluding}]
Let $\gamma\in \bar{\mathcal{E}}_0$, using the area formula (\cite[Theorem~3.2.3]{Fed69}) we obtain that
\[
L_0(\gamma)=\|\dot{\gamma}\|_\infty \ge \|\dot{\gamma}\|_1=\int\limits_{\RZ}|\dot{\gamma}|dt=\int\limits_{\gamma(\RZ)} \#(\gamma^{-1}(x))d\H^1(x).
\]
\underline{Claim}: $\#({\gamma}^{-1}(x))$ is even for $\H^1$-a.e. $x\! \in\! \gamma(\RZ)$ and in particular $\#({\gamma}^{-1}(x))\!\ge\! 2$.\\

Given a $1$-form $\alpha=\alpha_1 dx_1 + \alpha_2 dx_2$ on $\R^2$ of class $C^1$, by a parametric Gauss-Green theorem (see e.g., \cite[Theorem~1]{OrtelSchneider89}) we obtain, using the winding number hypothesis on $\gamma$, that
\begin{eqnarray*}
\int\limits_{\gamma} \alpha &:=&\int\limits_{0}^{1}\langle \alpha(\gamma(t)),\dot{\gamma}(t)\rangle dt\\
&=&\int\limits_{\R^2} \omega(\gamma, x)\left(\frac{\partial \alpha_2}{\partial x_1}-\frac{\partial\alpha_1}{\partial x_2}\right) dx\\ %=:\int\limits_{\R^2}\omega(\gamma,\alpha)d\alpha;
&=& 0.
\end{eqnarray*}
Now by the well known change of variable formula (see e.g. \cite[Section 5.4]{FonsecaGangbo}), we have that
\begin{equation}\label{Eq:deg=0}
0=\int\limits_{\gamma}\alpha = \int\limits_{\gamma(\RZ)}\langle \alpha(x),\tau(x)\rangle \mathtt{deg}(\gamma,x) d\H^1(x),
\end{equation}
where $\tau(x)$ is an orientation of the rectifable set $\gamma(\RZ)$ (see Remark~\ref{Rmk:Orientation}), and $\mathtt{deg}(\gamma,x)$ is defined for $\H^1$-a.e. $x\in \gamma(\RZ)$ as
\[
\mathtt{deg}(\gamma,x):=\sum\limits_{t\in\gamma^{-1}(x)} \pm 1,
\]
where $\pm 1$ is such that $\frac{\dot{\gamma}(t)}{|\dot{\gamma}(t)|}=\pm \tau(x)$.\\
Since $\alpha$ is an arbitrary $1$-form of class $C^1$, \eqref{Eq:deg=0} implies that $\mathtt{deg}(\gamma,x)=0$ for $\H^1$-a.e. $x\in \gamma(\RZ)$; and it follows that for such $x$, $\#({\gamma}^{-1}(x))$ is even.

The converse statement follows from Theorem 4.4 in \cite{AlbertiOttolini17}.
\end{proof}

\begin{comment}
\begin{theorem}\label{Cor:Concluding}
Given a set $K_0$ which is a minimizer for \Cref{Prob:Plateau}, there exists $\gamma_0\in \bar{\mathcal{E}}_0$ such that 
\begin{equation}\label{Eq:norm-length Equiv.}
L_0(\gamma_0) =  2\H^1(K_0);
\end{equation}
and $\gamma_0$ is a minimizer for the functional $L_0$. On the other hand, given a minimizer $\gamma\in \bar{\mathcal{E}}_0$ of $L_0$, it's image $K:=\gamma(\RZ)$ is a solution to \Cref{Prob:Plateau}.
\end{theorem}
\end{comment}

We denote by $\mathcal{K}$ the class of compact connected subsets $K\subset\R^2$.

\begin{proof}[Proof of \Cref{Cor:Concluding}] The proof follows by direct minimality arguments following \Cref{Thm:Concluding}.

Let $K_0\in \mathcal{K}$ be a solution to 
\Cref{Prob:Plateau} and let $\gamma_0\in \bar{\mathcal{E}}_0$ from \Cref{Thm:Concluding} such that $K_0=\gamma_0(\RZ)$ and
\[
L_0(\gamma_0)= 2\H^1(K_0).
\]
 
Now let $\gamma\in \bar{\mathcal{E}}_0$ be arbitrary, and let $K:=\gamma(\RZ)$; then by \Cref{Thm:Concluding} it holds that $K\in\mathcal{K}$ and
\[
L_0(\gamma)\ge 2\H^1(K).
\]
We thus have that
\[
L_0(\gamma_0)=2\H^1(K_0)\le 2\H^1(K)\le L_0(\gamma).
\]
Hence $\gamma_0$ is a minimizer of $L_0$.

For the converse statement: Let $\gamma_0\in \bar{\mathcal{E}}_0$ be a minimizer of $L_0$, let $K_0:=\gamma_0(\RZ)$ and  $K\in\mathcal{K}$ arbitrary. Let $\gamma\in \bar{\mathcal{E}}_0$ as in \Cref{Thm:Concluding} such that $K=\gamma(\RZ)$ and $L_0(\gamma)= 2\H^1(K)$. By the minimality of $\gamma_0$ we obtain that 
\[
2\H^1(K_0)\le L_0(\gamma_0)\le L_0(\gamma)= 2\H^1(K).
\]
And the conclusion follows.
%If you argue by contradiction, for the first result, you contradict minimality. The converse is trivial
\end{proof}

\subsection*{Acknowledgements.} %% To write if needed
I am most grateful to Giovanni Alberti for motivating problem and for many helpful discussions. Many thanks to Almut Burchard for valuable suggestions with improving presentation of original draft. 
The hospitality of the ICTP-East African Institute for Fundamental Research hosted by the University of Rwanda is most acknowledged.

\addtocontents{toc}{\protect\setcounter{tocdepth}{2}}

\bibliography{CPP.bib}{}
\bibliographystyle{plain}
\end{document}